\documentclass[12pt,reqno]{amsart}
\usepackage{amsfonts}
\usepackage{amsmath,amsthm,amsfonts,amssymb}
\newtheorem{lemma}{Lemma}[section]
\newtheorem{theorem}{Theorem}[section]
\newtheorem*{corollary}{Corollary}
\newtheorem*{Main Theorem}{Main Theorem}

\def\bl{\begin{lemma}}
\def\bt{\begin{theorem}}
\def\bbt{\begin{Main Theorem}}
\def\el{\end{lemma}}
\def\et{\end{theorem}}
\def\eet{\end{Main Theorem}}
\def\bp{\begin{proof}}
\def\ep{\end{proof}}
\def\bc{\begin{corollary}}
\def\ec{\end{corollary}}

\def\mb{\mathbb}

\def\a{\alpha}
\def\b{\beta}
\def\p{\partial}

\def\-{\setminus}

\def\vp{\varphi}

\def\lt{\left}
\def\rt{\right}
\def\+{\bigcup}
\def\.{\bigcap}

\title[]
{On the Poincar\'e-Lelong equation in $\mb{C}^n$}

\author []{Shaoyu Dai$^1$ and Yifei Pan$^2$}

\address{1 Department of Mathematics, Jinling Institute of Technology, Nanjing, 211169, China.}
\address{\it E-mail address: dymdsy@163.com}

\address{2 Department of Mathematical Sciences, Purdue University Fort Wayne, Fort Wayne, 46805-1499, USA.}
\address{\it E-mail address: pan@pfw.edu}

\begin{document}

\begin{abstract}
In this paper, we prove the existence of (global) solutions of the Poincar\'e-Lelong equation $\p\overline{\p}u=f$, where $f$ is a $d$-closed $(1,1)$ form and is in the weighted Hilbert space with Gaussian measure, i.e., $L^2_{(1,1)}(\mb{C}^n,e^{-|z|^2})$. The novelty of this paper is to apply a weighted $L^2$ version of Poincar\'e Lemma for $2$-forms, and then apply H\"{o}rmander's $L^2$ solutions for Cauchy-Riemann equations. In the both cases, the same weight $e^{-|z|^2}$ is used.
\end{abstract}

\maketitle

\section{Introduction}

In this paper, we will study the Poincar\'e-Lelong equation and prove the existence of (global) solutions in a weighted Hilbert space in $\mb{C}^n$, where $n\geq1$. More precisely, we prove the following theorem.

\bbt
For each $(1,1)$ form $f$ in the weighted Hilbert space $ L^2_{(1,1)}(\mb{C}^n,e^{-|z|^2})$ with $\p f=\overline{\partial}f=0$, there exists a solution $u$ in $L^2(\mb{C}^n,e^{-|z|^2})$ solving the Poincar\'e-Lelong equation
$$\partial\overline{\partial}u=f$$ in $\mb{C}^n$, in the sense of distributions,
with the norm estimate $$\int_{\mb{C}^n} |u|^2e^{-|z|^2}\leq
2\int_{\mb{C}^n} |f|^2e^{-|z|^2}.$$
\eet

Since $L^2_{(1,1)}(\mb{C}^n)$ or
$L^\infty_{(1,1)}(\mb{C}^n)\subset L^2_{(1,1)}(\mb{C}^n,e^{-|z|^2})$, we have the following corollary.

\bc
For each $(1,1)$ form $f\in L^2_{(1,1)}(\mb{C}^n)$ or $f\in L^\infty_{(1,1)}(\mb{C}^n)$ with $\p f=\overline{\partial}f=0$, there exists a solution $u$ in $L^2_{loc}(\mb{C}^n)$ solving the Poincar\'e-Lelong equation
$$\partial\overline{\partial}u=f$$ in $\mb{C}^n$, in the sense of distributions. In particular, if $f$ is further a positive $(1,1)$ form, then the solution must be plurisubharmonic.
\ec

Because of the main theorem, we have proved, in particular, that as far as the solvability is concerned, the Poincar\'e-Lelong equation can be solved globally for any $d$-closed $(1,1)$ form in the union of Hilbert spaces:  $\bigcup_{\lambda>0}L^2_{(1,1)}(\mb{C}^n,e^{-\lambda|z|^2})$.

P. Lelong \cite{1} studied the equation $\p\overline{\p}u=f$ in connection with questions on entire functions, and showed, unexpectedly, that with suitable restrictions on the growth of $f$,
the equation could be reduced to solving the more familiar equation $\frac{1}{4}\Delta u=\mbox{trace}(f)$ (Poisson equation). Mok, Siu and Yau \cite{2} studied the equation on a complete K\"{a}hler manifold and obtained important applications to questions on when a (noncompact) K\"{a}hler manifold is biholomorphicly equivalent to $\mb{C}^n$. Recently, Chen \cite{3} obtained solutions of the equation when $f$ is assumed to be a smooth $(1,1)$ $d$-closed form with compact support in $\mb{C}^n$, and he applied the result to prove a version of Hartog's extension theorem for pluriharmonic functions.

It was Berndtsson \cite{5}, who first studied the $d$-equation for $1$-forms and pointed out that
the H\"{o}rmander's $L^2$ method could be used for the $d$-equation in convex domains and with a convex weight function. Since our proof of the main theorem depends significantly on a weighted $L^2$ version of Poincar\'e Lemma (below), and the classical Poincar\'e Lemma would not provide a $L^2$ estimates for $d$-equation, we decide to include a detailed proof of Poincar\'e Lemma despite of the fact that we only use the weight $e^{-|x|^2}$ in solving the Poincar\'e-Lelong equation. In addition, the proof will provide a specific constant that we shall use in the main theorem.

\bigskip
\noindent\textbf{Poincar\'e Lemma.} (A weighted $L^2$ version for $p+1$-forms)
\textit{
Let $n\geq1$ be an integer. Let $\vp$ be a strictly convex smooth function in $\mathbb{R}^n$ and there exists a constant $c>0$ such that
$$\sum^n_{j,k=1}
\frac{\p^2\vp}{\p x_j\p x_k}\omega_{j}\omega_{k}\geq c|\omega|^2$$ for all $\omega=(\omega_1,\cdots,\omega_n)\in\mathbb{R}^n$. Let $p$ be an integer with $0\leq p\leq n-1$.
Then, for each $f$, a $d$-closed $p+1$-form in the weighted Hilbert space $ L^2_{p+1}(\mb{R}^n,e^{-\vp})$, there exists a solution $u$ in $L^2_{p}(\mb{R}^n,e^{-\vp})$ solving equation
$$du=f$$ in $\mb{R}^n$, in the sense of distributions,
with the norm estimate $$\int_{\mb{R}^n} |u|^2e^{-\vp}\leq
\frac{1}{c(p+1)}\int_{\mb{R}^n}|f|^2e^{-\vp}.$$
}

By the same idea of this paper, we could solve Poincar\'e-Lelong equation over any convex domain in $\mathbb{C}^n$ with an appropriate convex weight, which we will return to in another paper in a near future.

This paper is rather self-contained and much of its length is devoted to the proof of Poincar\'e Lemma. In the final section of the paper, we will explain why we need Poincar\'e Lemma in proving the main theorem.

\bigskip
\noindent\textbf{Acknowledgements.} The second named author's research is in part supported by a grant from 2019 Pippert Science Research Scholar at Purdue University Fort Wayne and he is thankful for the Pippert  family's generosity and continuing support for science research.

\section{Preliminary}

Here, let $n\geq1$ and $p$ be an nonnegative integer. For multiindex $I=(i_1,\cdots,i_p)$, where $i_1,\cdots,i_p$ are integers between $1$ and $n$, define $|I|=p$ and $dx^I=dx_{i_1}\wedge\cdots\wedge dx_{i_p}$.
In general, a $p$-form $f$ is a formal combination
$$f=\sum'_{|I|=p}f_Idx^I,$$
where ${\sum}'$ implies that the summation is performed only over strictly increasing multi-indices and
$f_I:\mb{R}^n\rightarrow\mb{R}$ is a function for all $f_I$.
For $p$-forms $f$ and $g$, we denote by $f\cdot g$ their pointwise scalar product, i.e.,
$$f\cdot g=\sum'_{|I|=p}f_Ig_I.$$
Let $\vp$ be a smooth and nonnegative function on $\mb{R}^n$.
We consider the weighted Hilbert space for $p$-forms
$$L^2_p(\mb{R}^n,e^{-\vp})
=\{f=\sum'_{|I|=p}f_Idx^I\mid f_I\in L^2_{loc}(\mb{R}^n); \int_{\mb{R}^n}|f|^2 e^{-\vp}<+\infty\},$$
where $|f|^2=f\cdot f$. Note that
$$L^2_0(\mb{R}^n,e^{-\vp})=L^2(\mb{R}^n,e^{-\vp})
=\{f\mid f\in L^2_{loc}(\mb{R}^n); \int_{\mb{R}^n}|f|^2 e^{-\vp}<+\infty\}.$$
We denote
the weighted inner product for $f,g\in L^2_p(\mb{R}^n,e^{-\vp})$ by
$$\langle f,g\rangle_{L^2_p(\mb{R}^n,e^{-\vp})}=\int_{\mb{R}^n}
f\cdot g e^{-\vp},$$
and the weighted norm of $f\in L^2_p(\mb{R}^n,e^{-\vp})$ by
$\|f\|_{L^2_p(\mb{R}^n,e^{-\vp})}=\sqrt{\langle f,f\rangle}_{L^2_p(\mb{R}^n,e^{-\vp})}.$

In Section 2, 3 and 4, since we only deal with $\mb{R}^n$, we will simply write $L^2_p(e^{-\vp})$ etc, for the weighted $L^2$-spaces.
Let $\mathcal{D}_p$ denote the set of $p$-forms whose coefficients are smooth functions with compact support in $\mb{R}^n$.

For each $p$-form $u={\sum'_{|I|=p}}u_Idx^I\in L^2_p(e^{-\vp})$, in the sense of distributions, the differential $du$ is that: $du=0$ when $p=n$, and when $p<n$,
\begin{align*}
du&=\sum'_{|I|=p}\sum^n_{j=1}\frac{\p u_I}{\p x_j}dx_j\wedge dx^I\\
&=\sum'_{|I|=p}\sum_{j\notin I}\frac{\p u_I}{\p x_j}dx_j\wedge dx^I\\
&=\sum'_{|I|=p}\sum_{j\notin I}\frac{\p u_I}{\p x_j}\epsilon^{jI}_{(jI)'}dx^{(jI)'}\\
&=\sum'_{|M|=p+1}\lt(\sum_{j\in M}\frac{\p u_{M^j}}{\p x_j}\epsilon^{jM^j}_{M}\rt)dx^M,
\end{align*}
where $jI=(j,i_1,\cdots,i_p)$, $(jI)'$ is the permutation of $jI$ such that $(jI)'$ is a strictly increasing multiindex, $\epsilon^{jI}_{(jI)'}$ is the signature of the permutation (for example, the signature is $-1$ if only two indices are interchanged), and $M^j$ is the increasing multiindex with $j$ removed from $M$.
For $p$-form $u\in L^2_p(e^{-\vp})$ and $p+1$-form $f\in L^2_{p+1}(e^{-\vp})$, we say that $f$ is the differential $du$ (in the sense of distributions), written $du=f$, provided
\begin{align*}
\int_{\mb{R}^n} du\cdot\a=\int_{\mb{R}^n} f\cdot\a
\end{align*}
for all test forms $\a={\sum'_{|J|=p+1}}\a_Jdx^J$ in $\mathcal{D}_{p+1}$.

Obviously, the operator
$d$ is well defined on $\mathcal{D}_p$:
$$d: \mathcal{D}_p\rightarrow \mathcal{D}_{p+1}.$$
We now extend the definition of the operator $d$ by allowing it to act on any $u\in L^2_p(e^{-\vp})$ such that $du$ (computed in the sense of distributions) lies in $L^2_{p+1}(e^{-\vp})$. This way we get a closed, densely defined operator
$$T: L^2_p(e^{-\vp})\rightarrow L^2_{p+1}(e^{-\vp}),$$
where the domain of $T$ is
$$Dom(T)=\{u\in L^2_p(e^{-\vp})\mid du\in L^2_{p+1}(e^{-\vp})\}.$$

Now we consider the Hilbert space adjoint of $T$:
$$T^*: L^2_{p+1}(e^{-\vp})\rightarrow L^2_p(e^{-\vp}).$$
Let $Dom(T^*)$ be the domain of $T^*$. Let $\a\in L^2_{p+1}(e^{-\vp})$. By functional analysis, we say that $\a\in Dom(T^*)$ if there exists a constant $c=c(\a)>0$ such that $$|\langle Tu,\a\rangle_{L^2_{p+1}(e^{-\vp})}|\leq c\|u\|_{L^2_p(e^{-\vp})}$$
for all $u\in Dom(T)$. This definition is equivalent to that  $\a\in Dom(T^*)$ if and only if there exists $v\in L^2_p(e^{-\vp})$ such that
$$\langle u,v\rangle_{L^2_p(e^{-\vp})}=\langle Tu,\a\rangle_{L^2_{p+1}(e^{-\vp})}$$
for all $u\in Dom(T)$. Note that $v$ is unique. We set $v=T^*\a$. Then
$T^*: Dom(T^*)\rightarrow L^2_p(e^{-\vp})$ is a linear operator and satisfies
\begin{align}\label{1}
\langle u,T^*\a\rangle_{L^2_p(e^{-\vp})}=\langle Tu,\a\rangle_{L^2_{p+1}(e^{-\vp})}
\end{align}
for all $u\in Dom(T)$, $\a\in Dom(T^*)$. It is well-known that $T^*$ is again a closed, densely defined operator.

In order to compute $T^*$, we first
computer $T^*_{formal}$, the
formal adjoint of $T$, which is defined using only test forms, i.e., we demand
\begin{align}\label{10000}
\langle Tu,\a\rangle_{L^2_{p+1}(e^{-\vp})}=\langle u,T^*_{formal}\a\rangle_{L^2_p(e^{-\vp})}
\end{align}
for $u\in Dom(T)$ and $\a\in \mathcal{D}_{p+1}$.
Note that for $u={\sum'_{|I|=p}}u_Idx^I$,
\begin{align*}
Tu=du=\sum'_{|I|=p}\sum^n_{j=1}\frac{\p u_I}{\p x_j}dx_j\wedge dx^I.
\end{align*}
For $\a=\sum'_{|J|=p+1}\a_Jdx^J$, if $J_1$ is a permutation of $J$, we write $\a_{J_1}=\epsilon^{J_1}_J\a_J$, where $\epsilon$ is the signature of the permutation. In particular, a term $\a_{jK}=0$ if $j\in K$. Then by
integration by parts, the left side of (\ref{10000}) is given by
\begin{align*}
\langle Tu,\a\rangle_{L^2_{p+1}(e^{-\vp})}&=\int_{\mb{R}^n}
du\cdot\a e^{-\vp}\\
&=\int_{\mb{R}^n}\sum'_{|I|=p}\sum^n_{j=1}\frac{\p u_I}{\p x_j}\a_{jI}e^{-\vp}\\
&=-\int_{\mb{R}^n}\sum'_{|I|=p}\sum^n_{j=1}u_I\frac{\p(\a_{jI}e^{-\vp})}{\p x_j}\\
&=\int_{\mb{R}^n}\lt(\sum'_{|I|=p}u_I\lt(-e^{\vp}
\sum^n_{j=1}\frac{\p(\a_{jI}e^{-\vp})}{\p x_j}\rt)\rt)e^{-\vp}\\
&=\int_{\mb{R}^n}\lt({\sum_{|I|=p}}'u_IA_I\rt)e^{-\vp}
,
\end{align*}
where
\begin{align}\label{5}
A_I=-e^{\vp}
\sum^n_{j=1}\frac{\p(\a_{jI}e^{-\vp})}{\p x_j}.
\end{align}
For example, if $p=1$, then
$$A_I=A_i=
-e^{\vp}\lt(\sum_{1\leq j<i}\frac{\p(\a_{ji}e^{-\vp})}{\p x_j}-\sum_{i<j\leq n}\frac{\p(\a_{ij}e^{-\vp})}{\p x_j}\rt).$$
Clearly, $A_I$ is a smooth function with compact support in $\mb{R}^n$, So ${\sum'_{|I|=p}}A_Idx^I\in\mathcal{D}_p\subset L^2_p(e^{-\vp})$.
Thus, the formal adjoint is
\begin{align}\label{6}
T^*_{formal}\a=\sum'_{|I|=p}A_Idx^I,
\end{align}
where $A_I$ is as (\ref{5}). This implies that $\mathcal{D}_{p+1}\subset Dom(T^*)$.

In the sense of distributions, the formal adjoint $T^*_{formal}\a$ is actually well-defined for $\forall\a\in L^2_{p+1}(e^{-\vp})$  as $\vp$ is smooth.
We claim that
\begin{align}\label{3}
T^*\a=T^*_{formal}\a\ \ \mbox{for}\ \ \forall\a\in Dom(T^*).
\end{align}
Indeed, if $\a\in Dom(T^*)$, then by (\ref{1}) and
$\mathcal{D}_p\in Dom(T)$, we have for $\forall u\in \mathcal{D}_p$,
$$\int_{\mb{R}^n}
u\cdot T^*\a e^{-\vp}=\langle u,T^*\a\rangle_{L^2_p(e^{-\vp})}=\langle Tu,\a\rangle_{L^2_{p+1}(e^{-\vp})}=\int_{\mb{R}^n}
u\cdot T^*_{formal}\a e^{-\vp}.$$
Then (\ref{3}) is hold.

\section{Approximation}
In this section, let $0\leq p\leq n-1$. We will prove that the set of smooth $p+1$-forms with compact support is dense in $Dom(T^*)\cap Dom(S)$ in the graph norm. The argument is standard, and for completeness of the paper, we  include the detailed proofs. At the same time, we follow the arguments of Fornaess's lecture notes \cite{6} closely in the case of the H\"{o}rmander's $L^2$ $\overline{\p}$ estimates.

Consider the spaces $G_1=L^2_p(\mb{R}^n,e^{-\vp}),$ $G_2=L^2_{p+1}(\mb{R}^n,e^{-\vp}),$ and $G_3=L^2_{p+2}(\mb{R}^n,e^{-\vp})$
with $d$ operators $T: G_1\rightarrow G_2$ and $S: G_2\rightarrow G_3$. Let $Dom(S)$ be the domain of $S$.

Let $\lambda_\nu$ be a sequence of smooth functions on $[0,+\infty)$ such that $\lambda_\nu(t)=1$ for $0\leq t\leq\nu$, $0<\lambda_\nu(t)<1$ for
$\nu<t<\nu+1$, $\lambda_\nu(t)=0$ for $t\geq\nu+1$ and $|\lambda'_\nu(t)|\leq2$ for $\forall t\in[0,+\infty)$.
Let $$\eta_\nu(x)=\lambda_\nu(|x|)\ \ \mbox{for}\ \ \forall x\in\mb{R}^n.$$
Obviously, $\eta_\nu(x)$ is a real-valued smooth function on $\mb{R}^n$ and $\eta_\nu(x)=1$ for $|x|\leq\nu$, $0<\eta_\nu(x)<1$ for $\nu<|x|<\nu+1$, $\eta_\nu(x)=0$ for $|x|\geq\nu+1$.
Then we have $$|d\eta_\nu(x)|^2\leq4n\ \ \mbox{for}\ \ \forall x\in\mb{R}^n,$$ since $d\eta_\nu(x)=0$ for $|x|\leq\nu$ and
$$|d\eta_\nu(x)|^2=\sum^n_{k=1}\lt(\frac{\p\lambda_\nu(|x|)}{\p x_k}\rt)^2
=\sum^n_{k=1}(\lambda'_\nu(|x|))^2\lt(\frac{x_k}{|x|}\rt)^2\leq4n$$
for $|x|>\nu$.

\bl\label{5.1}
Let $f\in Dom(S)$. Then the sequence $\eta_\nu f\rightarrow f$ in $G_2$. Moreover $\eta_\nu f\in Dom(S)$  and $S(\eta_\nu f)\rightarrow S(f)$ in $G_3$.
\el

\bp
The sequence $|\eta_\nu f|\leq|f|$ and $\eta_\nu f$ converges pointwise to $f$, so by the Lebesgue dominated convergence theorem, $\int_{\mb{R}^n}|\eta_\nu f-f|^2e^{-\vp}\rightarrow0$ as $\nu\rightarrow+\infty$.

Since
$$d(\eta_\nu f)=\eta_\nu df+d\eta_\nu\wedge f$$
in the sense of distributions and $df=Sf\in G_3$, so to show that $\eta_\nu f\in Dom(S)$ we need to show that $d\eta_\nu\wedge f\in G_3$. For
$f=\sum'_{|J|=p+1}f_Jdx^J\in G_2$, we have that: $d\eta_\nu\wedge f=0$ when $p+1=n$, and when $p+1<n$,
\begin{align*}
d\eta_\nu\wedge f&=\lt(\sum^n_{j=1}\frac{\p\eta_\nu}{\p x_j}dx_j\rt)\wedge f\\
&=\sum'_{|J|=p+1}\sum^n_{j=1}\frac{\p\eta_\nu}{\p x_j}f_J dx_j\wedge dx^J\\
&=\sum'_{|J|=p+1}\sum_{j\notin J}\frac{\p\eta_\nu}{\p x_j}f_J dx_j\wedge dx^J\\
&=\sum'_{|J|=p+1}\sum_{j\notin J}\frac{\p\eta_\nu}{\p x_j}f_J\epsilon^{jJ}_{(jJ)'}dx^{(jJ)'}\\
&=\sum'_{|M|=p+2}\sum_{j\in M}\frac{\p\eta_\nu}{\p x_j}f_{M^j}\epsilon^{jM^j}_{M}dx^{M}.
\end{align*}
Note that $|d\eta_\nu|^2\leq4n$. So we obtain that: $|d\eta_\nu\wedge f|^2=0$ when $p+1=n$, and when $p+1<n$,
\begin{align}
\lt|d\eta_\nu\wedge f\rt|^2&=\sum'_{|M|=p+2}\lt(\sum_{j\in M}\frac{\p\eta_\nu}{\p x_j}f_{M^j}\epsilon^{jM^j}_{M}\rt)^2\nonumber\\
&\leq\sum'_{|M|=p+2}\lt(\sum_{j\in M}\lt(\frac{\p\eta_\nu}{\p x_j}\rt)^2\rt)\lt(\sum_{j\in M}\lt(f_{M^j}\epsilon^{jM^j}_{M}\rt)^2\rt)\nonumber\\
&\leq\sum'_{|M|=p+2}|d\eta_\nu|^2|f|^2\leq c|d\eta_\nu|^2|f|^2\nonumber\\
&\leq 4nc|f|^2,\label{10005}
\end{align}
where c is a constant.
Then
$$\int_{\mb{R}^n}\lt|d\eta_\nu\wedge f\rt|^2e^{-\vp}
\leq 4nc\int_{\mb{R}^n}|f|^2e^{-\vp}<+\infty,$$
which means that $d\eta_\nu\wedge f\in G_3$. So we have $\eta_\nu f\in Dom(S)$.

Note that
$$S(\eta_\nu f)-Sf=S(\eta_\nu f)-\eta_\nu Sf+\eta_\nu Sf-Sf,$$
$\eta_\nu S(f)\rightarrow S(f)$ in $G_3$ and $S(\eta_\nu f)-\eta_\nu Sf=d\eta_\nu\wedge f\rightarrow0$ in $G_3$. Then $S(\eta_\nu f)\rightarrow S(f)$ in $G_3$.
\ep

\bl\label{5.2}
Let $f\in Dom(T^*)$. Then the sequence $\eta_\nu f\rightarrow f$ in $G_2$. Moreover $\eta_\nu f\in Dom(T^*)$  and $T^*(\eta_\nu f)\rightarrow T^*(f)$ in $G_1$.
\el

\bp
By the Lebesgue dominated convergence theorem, $\eta_\nu f\rightarrow f$ in $G_2$.

For $u\in Dom(T)$, we have
\begin{align*}
\langle Tu, \eta_\nu f\rangle_{G_2}
&=\langle \eta_\nu Tu, f\rangle_{G_2}\\
&=\langle T(\eta_\nu u)-d\eta_\nu\wedge u, f\rangle_{G_2}\\
&=\langle T(\eta_\nu u), f\rangle_{G_2}-\langle d\eta_\nu\wedge u, f\rangle_{G_2}\\
&=\langle \eta_\nu u, T^*(f)\rangle_{G_1}-\langle d\eta_\nu\wedge u, f\rangle_{G_2}\\
&=\langle u, \eta_\nu T^*(f)\rangle_{G_1}-\langle d\eta_\nu\wedge u, f\rangle_{G_2}.
\end{align*}
Note that for $u=\sum'_{|I|=p}u_Idx^I\in G_1$, just like (\ref{10005}), we have
\begin{align*}
\lt|d\eta_\nu\wedge u\rt|^2&\leq c|d\eta_\nu|^2|u|^2
\leq 4nc|u|^2,
\end{align*}
where c is a constant.
So
\begin{align*}
\lt\|d\eta_\nu\wedge u\rt\|_{G_2}^2=\int_{\mb{R}^n}\lt|d\eta_\nu\wedge u\rt|^2e^{-\vp}
\leq4nc\int_{\mb{R}^n}\lt|u\rt|^2e^{-\vp}
=4nc\lt\|u\rt\|_{G_1}^2.
\end{align*}
Then
\begin{align*}
\lt|\langle Tu, \eta_\nu f\rangle_{G_2}\rt|&\leq\|u\|_{G_1}\|\eta_\nu T^*(f)\|_{G_1}+\lt\|d\eta_\nu\wedge u\rt\|_{G_2}
\|f\|_{G_2}\\
&\leq\lt(\|T^*(f)\|_{G_1}+\sqrt{4nc}\|f\|_{G_2}\rt)\|u\|_{G_1}.
\end{align*}
By the definition of $Dom(T^*)$, we have $\eta_\nu f\in Dom(T^*)$.

Since $$T^*(\eta_\nu f)-T^*f=T^*(\eta_\nu f)-\eta_\nu T^*f+\eta_\nu T^*f-T^*f$$
and $\eta_\nu T^*f\rightarrow T^*f$ in $G_1$,
it is sufficient to prove $T^*(\eta_\nu f)-\eta_\nu T^*f\rightarrow0$ in $G_1$.
For $u\in Dom(T)$, we have
\begin{align*}
\langle u,T^*(\eta_\nu f)-\eta_\nu T^*f\rangle_{G_1}
&=\langle Tu, \eta_\nu f\rangle_{G_2}-\langle u, \eta_\nu T^*f\rangle_{G_1}\\
&=\langle \eta_\nu Tu, f\rangle_{G_2}-\langle u, \eta_\nu T^*f\rangle_{G_1}\\
&=\langle T(\eta_\nu u)-d\eta_\nu\wedge u, f\rangle_{G_2}-\langle u, \eta_\nu T^*f\rangle_{G_1}\\
&=\langle \eta_\nu u, T^*f\rangle_{G_2}-\langle d\eta_\nu\wedge u, f\rangle_{G_2}-\langle u, \eta_\nu T^*f\rangle_{G_1}\\
&=\langle u, \eta_\nu T^*f\rangle_{G_2}-\langle d\eta_\nu\wedge u, f\rangle_{G_2}-\langle u, \eta_\nu T^*f\rangle_{G_1}\\
&=-\langle d\eta_\nu\wedge u, f\rangle_{G_2}.
\end{align*}
Note that $\mathcal{D}_p\in Dom(T)$.
Then for $\forall u\in\mathcal{D}_p$,
$$\lt|\int_{\mb{R}^n}\lt(T^*(\eta_\nu f)-\eta_\nu T^*f\rt)\cdot ue^{-\vp}\rt|
\leq\int_{\mb{R}^n}|d\eta_\nu\wedge u||f|e^{-\vp}
\leq\int_{\mb{R}^n}\sqrt{c}|d\eta_\nu||u||f|e^{-\vp},$$
which means that
$$\lt|T^*(\eta_\nu f)-\eta_\nu T^*f\rt|e^{-\vp}\leq\sqrt{c}|d\eta_\nu||f|e^{-\vp}$$
almost everywhere in $\mb{R}^n$.
Thus, for almost everywhere in $\mb{R}^n$, $\lt|T^*(\eta_\nu f)-\eta_\nu T^*f\rt|\rightarrow0$ and
$$\lt|T^*(\eta_\nu f)-\eta_\nu T^*f\rt|^2e^{-\vp}\leq4nc|f|^2e^{-\vp}.$$
So it follows from the Lebesgue dominated convergence theorem that $T^*(\eta_\nu f)\rightarrow \eta_\nu T^*f$ in $G_1$.
\ep

Next we will study smoothing. The following are two well-known  smoothing theorems \cite{6}.

\bl
Let $\chi$ be a smooth function with compact support in $\mb{R}^n$ and $\int_{\mb{R}^n}\chi(x)dx=1$. Set $\chi_\varepsilon(x)=\frac{1}{\varepsilon^N}
\chi(\frac{x}{\varepsilon})$. If $g\in L^2(\mb{R}^n)$, then the convolution $g\ast\chi_\varepsilon$ satisfies
\begin{align*}
(g\ast\chi_\varepsilon)(x)=\int_{\mb{R}^n}g(y)\chi_\varepsilon
(x-y)dy=\int_{\mb{R}^n}g(x-y)\chi_\varepsilon
(y)dy=\int_{\mb{R}^n}g(x-\varepsilon y)\chi
(y)dy
\end{align*}
and is a smooth function such that $\|g\ast\chi_\varepsilon-g\|_{L^2}\rightarrow0$ when $\varepsilon\rightarrow0$. The support of $g\ast\chi_\varepsilon$ has no points at distance $>\varepsilon$ from the support of $g$ if the support of $\chi$ lies in the unit ball.
\el

\bl\label{5.4}
Let $f_1, \cdots, f_n\in L^1_{loc}(\mb{R}^n)$. Also suppose that the distribution $\sum^n_{j=1}\frac{\p f_j}{\p x_j}\in L^1_{loc}(\mb{R}^n)$. Then
$$\lt(\sum^n_{j=1}\frac{\p f_j}{\p x_j}\rt)\ast\chi_\varepsilon=\sum^n_{j=1}\frac{\p (f_j\ast\chi_\varepsilon)}{\p x_j}.$$
\el

\bl\label{why1}
Let $f\in Dom(S)$ have compact support in $\mb{R}^n$. Then $f\ast\chi_\varepsilon\rightarrow f$ in $G_2$, $f\ast\chi_\varepsilon\in Dom(S)$ and $S(f\ast\chi_\varepsilon)\rightarrow Sf$ in $G_3$.
\el

\bp
Let $f=\sum'_{|J|=p+1}f_Jdx^J$.
Since $f\in Dom(S)\subset G_2$ and $f$ has compact support, so by the smoothing theorem, $f_J\ast\chi_\varepsilon
\rightarrow f_J$ in $L^2(\mathbb{R}^n)$ for all $J$. Then
$f\ast\chi_\varepsilon\rightarrow f$ in $G_2$. Since $f\ast\chi_\varepsilon$ is smooth with compact support, $f\ast\chi_\varepsilon\in Dom(S)$. Furthermore,
note that: $Sf=df=0$ when $p+1=n$, and when $p+1<n$,
\begin{align*}
Sf=df=\sum'_{|M|=p+2}\lt(\sum_{j\in M}\frac{\p f_{M^j}}{\p x_j}\epsilon^{jM^j}_M\rt)dx^M.
\end{align*}
Then by Lemma \ref{5.4}, we have $(Sf)\ast\chi_\varepsilon=S(f\ast\chi_\varepsilon)$. By the smoothing theorem, $(Sf)\ast\chi_\varepsilon\rightarrow Sf$ in $G_3$. Therefore $S(f\ast\chi_\varepsilon)\rightarrow Sf$ in $G_3$.
\ep

Next we prove a similar lemma for $T^*$.

\bl\label{why2}
Let $f\in Dom(T^*)$ have compact support in $\mb{R}^n$. Then $f\ast\chi_\varepsilon\rightarrow f$ in $G_2$, $f\ast\chi_\varepsilon\in Dom(T^*)$ and $T^*(f\ast\chi_\varepsilon)\rightarrow T^*f$ in $G_1$.
\el

\bp
Let $f=\sum'_{|J|=p+1}f_Jdx^J$.
Since $f\in Dom(T^*)\subset G_2$ and $f$ has compact support, so by the smoothing theorem, $f_J\ast\chi_\varepsilon
\rightarrow f_J$ in $L^2(\mathbb{R}^n)$ for all $J$. Then
$f\ast\chi_\varepsilon\rightarrow f$ in $G_2$.
Since $f\ast\chi_\varepsilon$ is smooth with compact support, $f\ast\chi_\varepsilon\in Dom(T^*)$.
Furthermore,
note that
\begin{align*}
T^*f&=\sum'_{|I|=p}\lt(-e^{\vp}
\sum^n_{j=1}\frac{\p(f_{jI}e^{-\vp})}{\p x_j}\rt)dx^I=
\sum'_{|I|=p}\lt(-
\sum^n_{j=1}\frac{\p f_{jI}}{\p x_j}+\sum^n_{j=1}f_{jI}\frac{\p\vp}{\p x_j}\rt)dx^I\\
\end{align*}
Since
$f\ast\chi_\varepsilon=\sum'_{|J|=p+1}
\lt(f_J\ast\chi_\varepsilon\rt)dx^J$,
so by the above formula, we can write
$T^*(f\ast\chi_\varepsilon)=\sum'_{|I|=p}B_Idx^I$, where
\begin{align}\label{38}
B_I=-
\sum^n_{j=1}\frac{\p (f_{jI}\ast\chi_\varepsilon)}{\p x_j}+\sum^n_{j=1}(f_{jI}\ast\chi_\varepsilon)\frac{\p\vp}{\p x_j}.
\end{align}
Since $T^*f\in G_1$, so $\lt(-
\sum^n_{j=1}\frac{\p f_{jI}}{\p x_j}\rt)\in L^2_{loc}(\mb{R}^n)$. Then
the first term of the right side of (\ref{38}) can be written, using Lemma \ref{5.4}, as
$\lt(-
\sum^n_{j=1}\frac{\p f_{jI}}{\p x_j}\rt)\ast\chi_\varepsilon$ and converges to $\lt(-
\sum^n_{j=1}\frac{\p f_{jI}}{\p x_j}\rt)$ in $L^2(\mb{R}^n)$ by the smoothing Theorem. The second part converges to $\sum^n_{j=1}f_{jI}\frac{\p\vp}{\p x_j}$ in $L^2(\mb{R}^n)$ by the smoothing Theorem.
Therefore $T^*(f\ast\chi_\varepsilon)\rightarrow T^*f$ in $G_1$.
\ep

\bl\label{bijin}
Let $f\in Dom(T^*)\cap Dom(S)$. Then there exists a sequence $\{f_n\}\subset\mathcal{D}_{p+1}$ such that $f_n\in Dom(T^*)\cap Dom(S)$, $f_n\rightarrow f$ in $G_{2}$, $T^*f_n\rightarrow T^*f$ in $G_1$ and $Sf_n\rightarrow Sf$ in $G_3$.
\el

\bp
Let $\delta>0$. Using Lemma \ref{5.1} for $S$ and Lemma \ref{5.2} for $T^*$, we can let $\nu_0$ be large enough that
$$\|\eta_{\nu_0}f-f\|_{G_2},\|T^*(\eta_{\nu_0}f)-T^*f\|_{G_1},
\|S(\eta_{\nu_0}f)-Sf\|_{G_3}<\frac{\delta}{2}$$
and $\eta_{\nu_0}f\in Dom(T^*)\cap Dom(S)$. Then by Lemma \ref{why1} and \ref{why2}, we have
for $\varepsilon>0$ small enough, $\hat{f}=(\eta_{\nu_0}f\ast\chi_\varepsilon)$ is in $Dom(T^*)\cap Dom(S)$ and
$$\|\hat{f}-\eta_{\nu_0}f\|_{G_2},
\|T^*\hat{f}-T^*(\eta_{\nu_0}f)\|_{G_1},
\|S\hat{f}-S(\eta_{\nu_0}f)\|_{G_3}<\frac{\delta}{2}.$$
Thus,
$$\|\hat{f}-f\|_{G_2},\|T^*\hat{f}-T^*f\|_{G_1},
\|S\hat{f}-Sf\|_{G_3}<\delta.$$
\ep

\section{Proof of Poincar\'e Lemma}
In this section, we will give the proof of Poincar\'e Lemma. We start with some lemmas. Let $0\leq p\leq n-1$.

\bl\label{lemmaifif}
For each $f\in L_{p+1}^2(e^{-\vp})$, there exists a solution $u\in L_{p}^2(e^{-\vp})$ solving the equation
$$du=f$$ in $\mb{R}^n$, in the sense of distributions
with the norm estimate
$$\|u\|^2_{L_{p}^2(e^{-\vp})}\leq c$$
if and only if
$$|\langle f,\a\rangle_{L_{p+1}^2(e^{-\vp})}|^2\leq c\lt\|T^*\a\rt\|^2_{L_{p}^2(e^{-\vp})}, \ \ \forall\a\in \mathcal{D}_{p+1},$$
where $c$ is a constant.
\el

\bp
(Necessity) Note that $du=f\in L_{p+1}^2(e^{-\vp})$. Then we have $du=Tu$. For $\forall\a\in \mathcal{D}_{p+1}$, from the definition of $T^*$ and Cauchy-Schwarz inequality, we have
$$
|\langle f,\a\rangle_{L_{p+1}^2(e^{-\vp})}|^2=|\langle Tu,\a\rangle_{L_{p+1}^2(e^{-\vp})}|^2=\lt|\lt\langle u,T^*\a\rt\rangle_{L_{p}^2(e^{-\vp})}\rt|^2
\leq\|u\|^2_{L_{p}^2(e^{-\vp})}
\lt\|T^*\a\rt\|^2_{L_{p}^2(e^{-\vp})}.
$$
Note that $\|u\|^2_{L_{p}^2(e^{-\vp})}\leq c$. Then $|\langle f,\a\rangle_{L_{p+1}^2(e^{-\vp})}|^2\leq c\lt\|T^*\a\rt\|^2_{L_{p}^2(e^{-\vp})}$.

(Sufficiency) Consider the subspace
$$E=\lt\{T^*\a\mid\a\in \mathcal{D}_{p+1}\rt\}\subset L_{p}^2(e^{-\vp}).$$
Define a linear functional $L_f: E\rightarrow\mb{R}$ by
$$L_f\lt(T^*\a\rt)=\langle f,\a\rangle_{L_{p+1}^2(e^{-\vp})}
.$$
Since
$$\lt|L_f\lt(T^*\a\rt)\rt|=\lt|\langle f,\a\rangle_{L_{p+1}^2(e^{-\vp})}\rt|
\leq\sqrt{c}\lt\|T^*\a\rt\|_{L_{p}^2(e^{-\vp})},$$
then $L_f$ is a bounded functional on $E$. So by Hahn-Banach's extension theorem, $L_f$ can be extended to a linear functional $\widetilde{L}_f$ on $L_{p}^2(e^{-\vp})$
such that
\begin{equation}\label{31}
\lt|\widetilde{L}_f(g)\rt|\leq\sqrt{c}\lt\|g\rt\|
_{L_{p}^2(e^{-\vp})}, \ \ \forall g\in L_{p}^2(e^{-\vp}).
\end{equation}
Using the Riesz representation theorem for $\widetilde{L}_f$, there exists a unique $u_0\in L_{p}^2(e^{-\vp})$ such that
\begin{equation}\label{32}
\widetilde{L}_f(g)=\langle u_0,g\rangle_{L_{p}^2(e^{-\vp})}, \ \ \forall g\in L_{p}^2(e^{-\vp}).
\end{equation}

Now we prove $du_0=f$. For $\forall\a\in \mathcal{D}_{p+1}$, apply $g=T^*\a$ in (\ref{32}). Then
$$\widetilde{L}_f\lt(T^*\a\rt)=\lt\langle u_0,T^*\a\rt\rangle_{L_{p}^2(e^{-\vp})}=\lt\langle Tu_0,\a\rt\rangle_{L_{p+1}^2(e^{-\vp})}.$$
Note that
$$\widetilde{L}_f\lt(T^*\a\rt)=L_f\lt(T^*\a\rt)=\langle f,\a\rangle_{L_{p+1}^2(e^{-\vp})}.$$
Therefore,
$$\lt\langle Tu_0,\a\rt\rangle_{L_{p+1}^2(e^{-\vp})}=\langle f,\a\rangle_{L_{p+1}^2(e^{-\vp})}, \ \ \forall\a\in \mathcal{D}_{p+1}.$$
Thus, $Tu_0=f$, i.e., $du_0=f$.

Next we give a bound for the norm of $u_0$. Let $g=u_0$ in (\ref{31}) and (\ref{32}). Then we have
$$\|u_0\|^2_{L_{p}^2(e^{-\vp})}=\lt|\langle u_0,u_0\rangle_{L_{p}^2(e^{-\vp})}\rt|=\lt|\widetilde{L}_f(u_0)\rt|
\leq\sqrt{c}\lt\|u_0\rt\|_{L_{p}^2(e^{-\vp})}.$$
Therefore, $\|u_0\|_{L_{p}^2(e^{-\vp})}^2\leq c$.

Let $u=u_0$. So there exists $u\in L_{p}^2(e^{-\vp})$ such that
$du=f$ with $\|u\|_{L_{p}^2(e^{-\vp})}^2\leq c$.
\ep

\bl\label{da}
Let $\a=\sum'_{|J|=p+1}\a_Jdx^J$. Then
$$|d\a|^2=\sum'_{|J|=p+1}\sum^n_{j=1}\lt|\frac{\p\a_J}{\p x_j}\rt|^2-\sum'_{|I|=p}\sum^n_{j,k=1}\frac{\p\a_{kI}}{\p x_j}\frac{\p\a_{jI}}{\p x_k}.$$
\el

\bp
Note that $$d\a=\sum'_{|J|=p+1}\sum^n_{j=1}\frac{\p\a_J}{\p x_j}dx_j\wedge dx^J.$$
We prove the lemma by two cases.

Case 1: $p+1=n.$
In this case $d\a=0$ for type reasons. Recall that $\a_{jK}=0$ if $j\in K$. Then for the second term on the right side of the formula,
\begin{align*}
-\sum'_{|I|=p}\sum^n_{j,k=1}\frac{\p\a_{kI}}{\p x_j}\frac{\p\a_{jI}}{\p x_k}
&=-\sum'_{|I|=n-1}\sum_{j,k\notin I}\frac{\p\a_{kI}}{\p x_j}\frac{\p\a_{jI}}{\p x_k}=-\sum'_{|I|=n-1}\sum_{j\notin I}\lt|\frac{\p\a_{jI}}{\p x_j}\rt|^2\\
&=-\sum'_{|I|=n-1}\sum_{j\notin I}\lt|\frac{\p\lt(\epsilon^{jI}_{(jI)'}\a_{(jI)'}\rt)}{\p x_j}\rt|^2=-\sum'_{|I|=n-1}\sum_{j\notin I}\lt|\frac{\p\a_{(jI)'}}{\p x_j}\rt|^2\\
&=-\sum'_{|J|=n}\sum^n_{j=1}\lt|\frac{\p\a_{J}}{\p x_j}\rt|^2
=-\sum'_{|J|=p+1}\sum^n_{j=1}\lt|\frac{\p\a_{J}}{\p x_j}\rt|^2,
\end{align*}
which is the same as the first term on the right side of the formula except for sign. Then the formula is proved.

Case 2: $p+1<n.$ We can write
\begin{align*}
d\a=\sum'_{|M|=p+2}\sum_{j\in M}\frac{\p\a_{M^j}}{\p x_j}dx_j\wedge dx^{M^j}=\sum'_{|M|=p+2}\lt(\sum_{j\in M}\frac{\p\a_{M^j}}{\p x_j}\epsilon^{jM^j}_M\rt)dx^M.
\end{align*}
So we obtain that
\begin{align*}
|d\a|^2&=\sum'_{|M|=p+2}\lt(\sum_{j\in M}\frac{\p\a_{M^j}}{\p x_j}\epsilon^{jM^j}_M\rt)^2\\
&=\sum'_{|M|=p+2}\sum_{j,k\in M}\frac{\p\a_{M^j}}{\p x_j}\frac{\p\a_{M^k}}{\p x_k}\epsilon^{jM^j}_M\epsilon^{kM^k}_M\\
&=\sum'_{|M|=p+2}\sum_{j\in M}\lt|\frac{\p\a_{M^j}}{\p x_j}\rt|^2
+\sum'_{|M|=p+2}\sum_{\substack{j,k\in M\\j\neq k}}\frac{\p\a_{M^j}}{\p x_j}\frac{\p\a_{M^k}}{\p x_k}\epsilon^{jM^j}_M\epsilon^{kM^k}_M\\
&=\sum'_{|J|=p+1}\sum_{j\notin J}\lt|\frac{\p\a_J}{\p x_j}\rt|^2
+\sum'_{|I|=p}\sum_{\substack{j,k\notin I\\j\neq k}}\frac{\p\a_{(kI)'}}{\p x_j}\frac{\p\a_{(jI)'}}{\p x_k}\epsilon^{j(kI)'}_{(j(kI)')'}\epsilon^{k(jI)'}_{(k(jI)')'}\\
&=\sum'_{|J|=p+1}\sum_{j\notin J}\lt|\frac{\p\a_J}{\p x_j}\rt|^2
+\sum'_{|I|=p}\sum_{\substack{j,k\notin I\\j\neq k}}\frac{\p\a_{kI}}{\p x_j}\frac{\p\a_{jI}}{\p x_k}\epsilon^{kI}_{(kI)'}\epsilon^{jI}_{(jI)'}
\epsilon^{j(kI)'}_{(j(kI)')'}\epsilon^{k(jI)'}_{(k(jI)')'}\\
&=\sum'_{|J|=p+1}\sum_{j\notin J}\lt|\frac{\p\a_J}{\p x_j}\rt|^2
-\sum'_{|I|=p}\sum_{\substack{j,k\notin I\\j\neq k}}\frac{\p\a_{kI}}{\p x_j}\frac{\p\a_{jI}}{\p x_k}.
\end{align*}
Note that
\begin{align*}
\sum'_{|J|=p+1}\sum_{j\in J}\lt|\frac{\p\a_J}{\p x_j}\rt|^2
=\sum'_{|I|=p}\sum_{j\notin I}\lt|\frac{\p\a_{(jI)'}}{\p x_j}\rt|^2
=\sum'_{|I|=p}\sum_{j\notin I}\lt|\frac{\epsilon^{jI}_{(jI)'}\p\a_{jI}}{\p x_j}\rt|^2
=\sum'_{|I|=p}\sum_{j\notin I}\lt|\frac{\p\a_{jI}}{\p x_j}\rt|^2.
\end{align*}
Then
\begin{align*}
|d\a|^2&=\lt(\sum'_{|J|=p+1}\sum_{j\notin J}\lt|\frac{\p\a_J}{\p x_j}\rt|^2
+\sum'_{|J|=p+1}\sum_{j\in J}\lt|\frac{\p\a_J}{\p x_j}\rt|^2\rt)\\
&\ \ \ \ \ \ \ \ \ \ \ \ -\lt(\sum'_{|I|=p}\sum_{j\notin I}\lt|\frac{\p\a_{jI}}{\p x_j}\rt|^2+\sum'_{|I|=p}\sum_{\substack{j,k\notin I\\j\neq k}}\frac{\p\a_{kI}}{\p x_j}\frac{\p\a_{jI}}{\p x_k}\rt)\\
&=\sum'_{|J|=p+1}\sum^n_{j=1}\lt|\frac{\p\a_J}{\p x_j}\rt|^2
-\sum'_{|I|=p}\sum_{j,k\notin I}\frac{\p\a_{kI}}{\p x_j}\frac{\p\a_{jI}}{\p x_k}\\
&=\sum'_{|J|=p+1}\sum^n_{j=1}\lt|\frac{\p\a_J}{\p x_j}\rt|^2
-\sum'_{|I|=p}\sum^n_{j,k=1}\frac{\p\a_{kI}}{\p x_j}\frac{\p\a_{jI}}{\p x_k}.
\end{align*}
\ep

\bl\label{zy}
Let $\a={\sum'_{|J|=p+1}}\a_Jdx^J\in \mathcal{D}_{p+1}$. Then
\begin{align}\label{10003}
\|T^*\a\|^2_{L^2_p(e^{-\vp})}+\|d\a\|^2_{L^2_{p+2}(e^{-\vp})}
=\int_{\mb{R}^n}\sum'_{|I|=p}\sum^n_{j,k=1}
\frac{\p^2\vp}{\p x_j\p x_k}\a_{jI}\a_{kI}e^{-\vp}+\int_{\mb{R}^n}
\sum'_{|J|=p+1}\sum^n_{j=1}\lt|\frac{\p\a_J}{\p x_j}\rt|^2 e^{-\vp}.
\end{align}
In particular, if there exists a constant $c>0$ such that
$$\sum^n_{j,k=1}
\frac{\p^2\vp}{\p x_j\p x_k}\omega_{j}\omega_{k}\geq c|\omega|^2$$ for all $\omega=(\omega_1,\cdots,\omega_n)\in\mathbb{R}^n$, then
\begin{align}\label{10004}
\|T^*\a\|^2_{L^2_p(e^{-\vp})}+\|d\a\|^2_{L^2_{p+2}(e^{-\vp})}
\geq c(p+1)\|\a\|^2_{L^2_{p+1}(e^{-\vp})}.
\end{align}
\el

\bp
We first prove (\ref{10003}). Consider the expression
\begin{align}\label{34}
Q=\|T^*\a\|^2_{L^2_p(e^{-\vp})}=\langle T^*\a,T^*\a\rangle_{L^2_p(e^{-\vp})}=\langle TT^*\a,\a\rangle_{L^2_p(e^{-\vp})}.
\end{align}
By (\ref{6}) and (\ref{3}), we have
$$T^*\a={\sum_{|I|=p}}'A_Idx^I,$$
where $A_I$ is as (\ref{5}).
Then
$$TT^*\a=d
\lt(\sum'_{|I|=p}A_Idx^I\rt)
=\sum'_{|I|=p}\sum^n_{k=1}\frac{\p A_I}{\p x_k}dx_k\wedge dx^I.
$$
Let $$\delta_j=e^\vp\frac{\p}{\p x_j}e^{-\vp}
=\frac{\p}{\p x_j}-\frac{\p\vp}{\p x_j}.$$
Then
$$A_I=-\sum^n_{j=1}\delta_j\a_{jI}.$$
Observe that
$$\frac{\p}{\p x_k}\delta_j=\frac{\p^2}
{\p x_j\p x_k}-\frac{\p\vp}{\p x_j}\frac{\p}
{\p x_k}-\frac{\p^2\vp}
{\p x_j\p x_k}$$
and
$$\delta_j\frac{\p}{\p x_k}=\frac{\p^2}
{\p x_j\p x_k}-\frac{\p\vp}{\p x_j}\frac{\p}
{\p x_k}.$$
We have
$$\frac{\p}{\p x_k}\delta_j=\delta_j\frac{\p}{\p x_k}-\vp_{jk}.$$
Here $\vp_{jk}=\frac{\p^2\vp}
{\p x_j\p x_k}.$
So for $1\leq k\leq n$, we have
\begin{align*}
\frac{\p A_I}{\p x_k}=-\sum^n_{j=1}\lt(\frac{\p}{\p x_k}\delta_j\rt)\a_{jI}
=\sum^n_{j=1}\lt(\vp_{jk}\a_{jI}-\delta_j\frac{\p\a_{jI}}{\p x_k}\rt)
\end{align*}
Then by (\ref{34}), we have
\begin{align}
Q&=\int_{\mb{R}^n}TT^*\a\cdot\a e^{-\vp}\nonumber\\
&=\int_{\mb{R}^n}\sum'_{|I|=p}\sum^n_{j=1}\frac{\p A_I}{\p x_k}\a_{kI}e^{-\vp}\nonumber\\
&=\int_{\mb{R}^n}\sum'_{|I|=p}\sum^n_{j,k=1}
\vp_{jk}\a_{jI}\a_{kI}e^{-\vp}
+\int_{\mb{R}^n}\sum'_{|I|=p}\sum^n_{j,k=1}(-1)
\lt(\delta_j\frac{\p\a_{jI}}{\p x_k}\rt)\a_{kI}e^{-\vp}
\nonumber\\
&=Q_1+Q_2.\label{10001}
\end{align}
Observe that
\begin{align*}
Q_2=\int_{\mb{R}^n}\sum'_{|I|=p}\sum^n_{j,k=1}(-1)\frac{\p}{\p x_j}
\lt(e^{-\vp}\frac{\p\a_{jI}}{\p x_k}\rt)\a_{kI}
=\int_{\mb{R}^n}\sum'_{|I|=p}\sum^n_{j,k=1}\frac{\p\a_{kI}}{\p x_j}
\frac{\p\a_{jI}}{\p x_k}e^{-\vp}.
\end{align*}
So by Lemma \ref{da}, we have
\begin{align}\label{10002}
Q_2=\int_{\mb{R}^n}\sum'_{|J|=p+1}\sum^n_{j=1}\lt|\frac{\p\a_J}{\p x_j}\rt|^2 e^{-\vp}-\|d\a\|^2_{L^2_{p+2}(e^{-\vp})}.
\end{align}
Then (\ref{10003}) is proved by (\ref{34}), (\ref{10001}) and (\ref{10002}).

Now we prove (\ref{10004}). Observe that
\begin{align*}
\sum'_{|I|=p}\sum^n_{j,k=1}
\vp_{jk}\a_{jI}\a_{kI}&\geq\sum'_{|I|=p}c\sum^n_{j=1}\lt|\a_{jI}\rt|^2
=c\sum'_{|I|=p}\sum_{j\notin I}\lt|\a_{jI}\rt|^2\\
&=c\sum'_{|I|=p}\sum_{j\notin I}\lt|\epsilon^{jI}_{(jI)'}\a_{(jI)'}\rt|^2=c\sum'_{|I|=p}\sum_{j\notin I}\lt|\a_{(jI)'}\rt|^2\\
&=c\sum'_{|J|=p+1}\sum_{j\in J}\lt|\a_{J}\rt|^2
=c\sum'_{|J|=p+1}(p+1)\lt|\a_{J}\rt|^2\\
&=c(p+1)|\a|^2.
\end{align*}
Then for the first term on the right side of (\ref{10003}),
\begin{align*}
\int_{\mb{R}^n}\sum'_{|I|=p}\sum^n_{j,k=1}
\vp_{jk}\a_{jI}\a_{kI}e^{-\vp}\geq\int_{\mb{R}^n}c(p+1)|\a|^2e^{-\vp}
=c(p+1)\|\a\|^2_{L^2_{p+1}(e^{-\vp})}.
\end{align*}
Note that the second term on the right side of (\ref{10003}) is always nonnegative. Then (\ref{10004}) is proved.
\ep

Now we give the proof of Poincar\'e Lemma.

\bp
Let $N=\{f\mid f\in L^2_{p+1}(e^{-\vp});df=0\}$, which is a closed subspace of $L^2_{p+1}(e^{-\vp})$.
For each $\a$ in $\mathcal{D}_{p+1}$, clearly $\a\in L^2_{p+1}(e^{-\vp})$, so we can decompose $\a=\a^1+\a^2$, where $\a^1$ lies in $N$ and $\a^2$ is orthogonal to $N$. This implies that $\a^2$ is orthogonal to any form $Tu$, since $Tu\in N$. So by the definition of $Dom(T^*)$, we see that $\a^2$ lies in the domain of $T^*$ and $T^*\a^2=0$. Since $\a$ lies in the domain of $T^*$, it follows that $T^*\a=T^*\a^1$.

Note that $\a^1\in Dom(T^*)\cap Dom(S)$. Then by Lemma \ref{bijin}, there exists a sequence $\{\a_\nu\}\subset\mathcal{D}_{p+1}$ such that $\a_\nu\in Dom(T^*)\cap Dom(S)$, $\a_\nu\rightarrow\a^1$ in $L_{p+1}^2(e^{-\vp})$, $T^*\a_\nu\rightarrow T^*\a^1$ in $L_{p}^2(e^{-\vp})$, and $S\a_\nu\rightarrow S\a^1$ in $L_{p+2}^2(e^{-\vp})$.

For $\a_\nu\in\mathcal{D}_{p+1}$, by Lemma \ref{zy}, we have
\begin{align*}
\|T^*\a_\nu\|^2_{L^2_p(e^{-\vp})}+\|S\a_\nu\|^2_{L^2_{p+2}(e^{-\vp})}\geq
c(p+1)\|\a_\nu\|^2_{L^2_{p+1}(e^{-\vp})}.
\end{align*}
Let $\nu\rightarrow+\infty$, so
\begin{align*}
\|T^*\a^1\|^2_{L^2_p(e^{-\vp})}+\|S\a^1\|^2_{L^2_{p+2}(e^{-\vp})}\geq
c(p+1)\|\a^1\|^2_{L^2_{p+1}(e^{-\vp})},
\end{align*}
which means that
\begin{align*}
\|T^*\a^1\|^2_{L^2_p(e^{-\vp})}\geq
c(p+1)\|\a^1\|^2_{L^2_{p+1}(e^{-\vp})}
\end{align*}
since $S\a^1=0$.

By Cauchy-Schwarz inequality, we have
\begin{align*}
\lt|\langle f,\a^1\rangle_{L_{p+1}^2(e^{-\vp})}\rt|^2
&\leq\lt\|f\rt\|^2_{L_{p+1}^2(e^{-\vp})}
\lt\|\a^1\rt\|^2_{L_{p+1}^2(e^{-\vp})}\\
&=\lt(\frac{1}{c(p+1)}\lt\|f\rt\|^2_{L_{p+1}^2(e^{-\vp})}\rt)
\lt(c(p+1)\lt\|\a^1\rt\|^2_{L_{p+1}^2(e^{-\vp})}\rt)\\
&\leq\lt(\frac{1}{c(p+1)}\lt\|f\rt\|^2_{L_{p+1}^2(e^{-\vp})}\rt)
\|T^*\a^1\|^2_{L^2_p(e^{-\vp})}.
\end{align*}
Let $\widetilde{c}=\frac{1}{c(p+1)}\lt\|f\rt\|^2_{L_{p+1}^2(e^{-\vp})}$. Then
\begin{align*}
\lt|\langle f,\a^1\rangle_{L_{p+1}^2(e^{-\vp})}\rt|^2
\leq \widetilde{c}
\|T^*\a^1\|^2_{L^2_p(e^{-\vp})}, \ \ \forall\a\in \mathcal{D}_{p+1}.
\end{align*}
Note that $f\in N$.
Thus,
$$\lt|\langle f,\a\rangle_{L_{p+1}^2(e^{-\vp})}\rt|^2=\lt|\langle f,\a^1\rangle_{L_{p+1}^2(e^{-\vp})}\rt|^2\leq \widetilde{c}\lt\|T^*\a^1\rt\|^2_{L_{p}^2(e^{-\vp})}
=\widetilde{c}\lt\|T^*\a\rt\|^2_{L_{p}^2(e^{-\vp})}.$$
By Lemma \ref{lemmaifif}, there exists a solution $u\in L_{p}^2(e^{-\vp})$ solving the equation
$$du=f$$ in $\mb{R}^n$, in the sense of distributions
with the norm estimate
$$\|u\|^2_{L_{p}^2(e^{-\vp})}\leq \widetilde{c},$$
i.e.,
$$\int_{\mb{R}^n}|u|^2e^{-\vp}\leq
\frac{1}{c(p+1)}\int_{\mb{R}^n}|f|^2e^{-\vp}.$$
The proof is complete.
\ep

\section{Proof of the main theorem}
Here, let $\vp$ be a smooth and nonnegative function on $\mb{C}^n$. We consider the weighted Hilbert space
$$L^2(\mb{C}^n,e^{-\vp})
=\{u:\mb{C}^n\rightarrow\mb{C}\mid u\in L^2_{loc}(\mb{C}^n); \int_{\mb{C}^n}|u|^2 e^{-\vp}<+\infty\}.$$
We denote
the weighted inner product for $u,v\in L^2(\mb{C}^n,e^{-\vp})$ by
$$\langle u,v\rangle_{L^2(\mb{C}^n,e^{-\vp})}=\int_{\mb{C}^n}u\overline{v} e^{-\vp},$$
and the weighted norm of $u\in L^2(\mb{C}^n,e^{-\vp})$ by
$\|u\|_{L^2(\mb{C}^n,e^{-\vp})}=\sqrt{\langle u,u\rangle}_{L^2(\mb{C}^n,e^{-\vp})}.$

In general, a $(1,1)$ form $f$ is a formal combination
$$f=\sum^n_{i,j=1}f_{i\overline{j}}dz_i\wedge d\overline{z}_j,$$
where $f_{i\overline{j}}:\mb{C}^n\rightarrow\mb{C}$ is a function for $1\leq i<j\leq n$. For $(1,1)$ forms $f$ and $g$, we denote by $f\cdot \overline{g}$ their pointwise scalar product, i.e.,
$$f\cdot \overline{g}=\sum^n_{i,j=1}f_{i\overline{j}}\overline{g}_{i\overline{j}}.$$
We also consider the weighted Hilbert space for $(1,1)$ forms
$$L^2_{(1,1)}(\mb{C}^n,e^{-\vp})
=\{f=\sum^n_{i,j=1}f_{i\overline{j}}dz_i\wedge d\overline{z}_j\mid f_{i\overline{j}}\in L^2_{loc}(\mb{C}^n); \int_{\mb{C}^n}|f|^2 e^{-\vp}<+\infty\},$$
where $|f|^2=f\cdot \overline{f}$. We denote
the weighted inner product for $f,g\in L^2_{(1,1)}(\mb{C}^n,e^{-\vp})$ by
$$\langle f,g\rangle_{L^2_{(1,1)}(\mb{C}^n,e^{-\vp})}=\int_{\mb{C}^n}
f\cdot \overline{g} e^{-\vp},$$
and the weighted norm of $f\in L^2_{(1,1)}(\mb{C}^n,e^{-\vp})$ by
$\|f\|_{L^2_{(1,1)}(\mb{C}^n,e^{-\vp})}=\sqrt{\langle f,f\rangle}_{L^2_{(1,1)}(\mb{C}^n,e^{-\vp})}.$

First we give two lemmas concerning about the conversion between complex forms and real forms.

\bl\label{wj1}
Let $f\in L^2_{(1,1)}(\mb{C}^n,e^{-\vp})$. Then
$f$ can be decomposed to $$f=f_1+\sqrt{-1}f_2,$$
where $f_1,f_2\in L^2_{2}(\mb{R}^{2n},e^{-\vp})$. Moreover,
\begin{align*}
|f_1|^2+|f_2|^2=4|f|^2.
\end{align*}
\el

\bp
Let
$
f=\sum^n_{i,j=1}f_{i\overline{j}}dz_i\wedge d\overline{z}_j,
$
where $f_{i\overline{j}}:\mb{C}^n\rightarrow\mb{C}$ is a function for $1\leq i,j\leq n$.
Let
$f_{i\overline{j}}=A_{i\overline{j}}+\sqrt{-1}B_{i\overline{j}},$
where $A_{i\overline{j}}:\mb{C}^n\rightarrow\mb{R}$ is a function for $1\leq i,j\leq n$ and $B_{i\overline{j}}:\mb{C}^n\rightarrow\mb{R}$ is a function for $1\leq i,j\leq n$.
Let
$z_i=x_i+\sqrt{-1}y_i$ for $1\leq i\leq n$.
Then
$
dz_i\wedge d\overline{z}_j=\lt(dx_i\wedge dx_j+dy_i\wedge dy_j\rt)-\sqrt{-1}\lt(dx_i\wedge dy_j+d x_j\wedge dy_i\rt)
$
and
\begin{align*}
f_{i\overline{j}}dz_i\wedge d\overline{z}_j
&=A_{i\overline{j}}dx_i\wedge dx_j+A_{i\overline{j}}dy_i\wedge dy_j+B_{i\overline{j}}dx_i\wedge dy_j+B_{i\overline{j}}d x_j\wedge dy_i\nonumber\\
&\ \ \ \ +\sqrt{-1}\lt(B_{i\overline{j}}dx_i\wedge dx_j+B_{i\overline{j}}dy_i\wedge dy_j-A_{i\overline{j}}dx_i\wedge dy_j-A_{i\overline{j}}d x_j\wedge dy_i\rt).
\end{align*}
Thus, we have
$
f=f_1+\sqrt{-1}f_2,
$
where
\begin{align*}
f_1=\sum_{1\leq i<j\leq n}\lt(A_{i\overline{j}}-A_{j\overline{i}}\rt)dx_i\wedge dx_j
+\sum_{1\leq i<j\leq n}\lt(A_{i\overline{j}}-A_{j\overline{i}}\rt)dy_i\wedge dy_j +\sum^n_{i,j=1}\lt(B_{i\overline{j}}
+B_{j\overline{i}}\rt)dx_i\wedge dy_j
\end{align*}
and
\begin{align*}
f_2=\sum_{1\leq i<j\leq n}\lt(B_{i\overline{j}}-B_{j\overline{i}}\rt)dx_i\wedge dx_j+
\sum_{1\leq i<j\leq n}\lt(B_{i\overline{j}}-B_{j\overline{i}}\rt)dy_i\wedge dy_j
-\sum^n_{i,j=1}\lt(A_{i\overline{j}}
+A_{j\overline{i}}\rt)dx_i\wedge dy_j
.
\end{align*}
Obviously, $f_1$ and $f_2$ are $2$-forms in $\mb{R}^{2n}$. We have $f_1,f_2\in L^2_{2}(\mb{R}^{2n},e^{-\vp})$, since $f\in L^2_{(1,1)}(\mb{C}^n,e^{-\vp})$.

Moreover,
\begin{align*}
|f_1|^2=2\sum_{1\leq i<j\leq n}\lt(A_{i\overline{j}}-A_{j\overline{i}}\rt)^2 +\sum^n_{i,j=1}\lt(B_{i\overline{j}}
+B_{j\overline{i}}\rt)^2
=\sum^n_{i,j=1}\lt(\lt(A_{i\overline{j}}-A_{j\overline{i}}\rt)^2+
\lt(B_{i\overline{j}}
+B_{j\overline{i}}\rt)^2\rt)
\end{align*}
and
\begin{align*}
|f_2|^2=2\sum_{1\leq i<j\leq n}\lt(B_{i\overline{j}}-B_{j\overline{i}}\rt)^2+
\sum^n_{i,j=1}\lt(A_{i\overline{j}}
+A_{j\overline{i}}\rt)^2=\sum^n_{i,j=1}\lt(
\lt(B_{i\overline{j}}-B_{j\overline{i}}\rt)^2
+\lt(A_{i\overline{j}}
+A_{j\overline{i}}\rt)^2\rt).
\end{align*}
Then
\begin{align*}
|f_1|^2+|f_2|^2
=2\sum^n_{i,j=1}\lt(A_{i\overline{j}}^2
+A_{j\overline{i}}^2+B_{i\overline{j}}^2
+B_{j\overline{i}}^2\rt)=4\sum^n_{i,j=1}\lt(A_{i\overline{j}}^2
+B_{i\overline{j}}^2
\rt)=4|f|^2.
\end{align*}
\ep

\bl\label{wj2}
Let $v\in L^2_{1}(\mb{R}^{2n},e^{-\vp})$. Then
$v$ can be decomposed to $$v=v^{1,0}+v^{0,1},$$
where $v^{1,0}\in L^2_{{1,0}}(\mb{C}^{n},e^{-\vp})$, $v^{0,1}\in L^2_{{0,1}}(\mb{C}^{n},e^{-\vp})$, $\overline{v^{1,0}}=v^{0,1}$ and $\overline{v^{0,1}}=v^{1,0}$.
Moreover,
$$\lt|v^{1,0}\rt|^2=\lt|v^{0,1}\rt|^2=\frac{1}{4}|v|^2.$$
\el

\bp
Let $v=\sum^{2n}_{j=1}v_jdx_j,$
where $v_j:\mb{R}^{2n}\rightarrow\mb{R}$ is a function for $1\leq j\leq2n$.
Let $z_j=x_{2j-1}+\sqrt{-1}x_{2j}.$
Then
$v=\sum^{n}_{j=1}\lt(v_{2j-1}dx_{2j-1}+v_{2j}dx_{2j}\rt)
=v^{1,0}+v^{0,1}$,
where
$$
v^{1,0}=\sum^{n}_{j=1}\lt(\frac{1}{2}v_{2j-1}
+\frac{1}{2\sqrt{-1}}v_{2j}\rt)dz_j
\ \ \ \ \mbox{and}
\ \ \ \
v^{0,1}=\sum^{n}_{j=1}\lt(\frac{1}{2}v_{2j-1}
-\frac{1}{2\sqrt{-1}}v_{2j}\rt)d\overline{z}_j.
$$
Obviously, $v^{1,0}$ is a $(1,0)$ form in $\mb{C}^{n}$, $v^{0,1}$ is a $(0,1)$ form in $\mb{C}^{n}$, $\overline{v^{1,0}}=v^{0,1}$ and $\overline{v^{0,1}}=v^{1,0}$.
Since $v\in L^2_{1}(\mb{R}^{2n},e^{-\vp})$, we have $v^{1,0}\in L^2_{{1,0}}(\mb{C}^{n},e^{-\vp})$ and $v^{0,1}\in L^2_{{0,1}}(\mb{C}^{n},e^{-\vp})$.

Moreover,
$$\lt|v^{1,0}\rt|^2=\sum^{n}_{j=1}
\lt(\lt(\frac{1}{2}v_{2j-1}\rt)^2+\lt(\frac{1}{2}v_{2j}\rt)^2\rt)
=\frac{1}{4}\sum^{2n}_{j=1}v_j^2=\frac{1}{4}|v|^2.$$
Similarly,
$$\lt|v^{0,1}\rt|^2=\frac{1}{4}|v|^2.$$
\ep

Now we give three more lemmas.

\bl\label{pu}
If $u\in L^2(\mb{C}^{n},e^{-\vp})$ and $\overline{\p}u\in L^2_{{0,1}}(\mb{C}^{n},e^{-\vp})$. Then $\p\overline{u}=\overline{\overline{\p}u}$, where $\overline{\p}u$ and $\p\overline{u}$ are in the sense of distributions.
\el

\bp
Let $\overline{\p}u=\sum^n_{j=1}v_jd\overline{z}_j\in L^2_{{0,1}}(\mb{C}^{n},e^{-\vp})$, where $v_j=\frac{\p u}{\p\overline{z}_j}$ in the sense of distributions. Then $\overline{\overline{\p}u}=\sum^n_{j=1}\overline{v_j}dz_j\in L^2_{{1,0}}(\mb{C}^{n},e^{-\vp})$.
For any $(1,0)$ test form $\a=\sum^n_{j=1}\a_jdz_j$, whose coefficients are smooth functions with compact support in $\mb{C}^{n}$, we have
\begin{align*}
\overline{\overline{\p}u}(\a)=\int_{\mb{C}^{n}}\sum^n_{j=1}
\overline{v_j}\a_j=\overline{\int_{\mb{C}^{n}}\sum^n_{j=1}
v_j\overline{\a_j}}=\overline{-\int_{\mb{C}^{n}}u\sum^n_{j=1}
\frac{\p\overline{\a_j}}{\p\overline{z}_j}}
=-\int_{\mb{C}^{n}}\overline{u}\sum^n_{j=1}\frac{\p\a_j}{\p z_j}
=\p\overline{u}(\a).
\end{align*}
Then $\p\overline{u}=\overline{\overline{\p}u}$.
\ep

\bl\label{huhuan}
If $u\in L^2(\mb{C}^{n},e^{-\vp})$. Then $\overline{\p}\p u=-\p\overline{\p}u$ in the sense of distributions.
\el

\bp
For $u\in L^2(\mb{C}^{n},e^{-\vp})$, we have
$$\overline{\p}\p u=\overline{\p}\lt(\sum^n_{i=1}\frac{\p u}{\p z_i}dz_i\rt)=
\sum^n_{i,j=1}\frac{\p^2u}{\p\overline{z}_j\p z_i}d\overline{z}_j\wedge dz_i$$
and
$$\p\overline{\p} u=\p\lt(\sum^n_{j=1}\frac{\p u}{\p \overline{z}_j}dz_j\rt)=
\sum^n_{i,j=1}\frac{\p^2u}{\p z_i\p\overline{z}_j}dz_i\wedge d\overline{z}_j=-\sum^n_{i,j=1}\frac{\p^2u}{\p z_i\p\overline{z}_j}d\overline{z}_j\wedge dz_i,$$
where $\frac{\p^2u}{\p\overline{z}_j\p z_i}$ and $\frac{\p^2u}{\p z_i\p\overline{z}_j}$ are in the sense of distributions.
Note that $\frac{\p^2u}{\p\overline{z}_j\p z_i}=\frac{\p^2u}{\p z_i\p\overline{z}_j}$. Then $\overline{\p}\p u=-\p\overline{\p}u$.
\ep

\noindent\textbf{Remark 5.1.} In the lemma, it is crucial that $\overline{\p}\p u$ and $\p\overline{\p}u$ are both forms. Otherwise, when $n=1$, $\overline{\p}\p u=\p\overline{\p}u$ if $\overline{\p}\p u=\p\overline{\p}u=\frac{\p^2u}{\p z\p \overline{z}}$ are as weak derivatives.

\bl\label{zuijiao}
Let $u\in L^2(\mb{C}^{n},e^{-\vp})$. If $\overline{\p}u\in L^2_{{0,1}}(\mb{C}^{n},e^{-\vp})$, then $\p\overline{\p}u=\p(\overline{\p}u)$ in the sense of distributions. If $\p u\in L^2_{{1,0}}(\mb{C}^{n},e^{-\vp})$, then $\overline{\p}\p u=\overline{\p}(\p u)$ in the sense of distributions.
\el

\bp
If $\overline{\p}u\in L^2_{{0,1}}(\mb{C}^{n},e^{-\vp})$, then for any $(1,1)$ test form $\a=\sum^n_{i,j=1}\a_{i\overline{j}}dz_i\wedge d\overline{z}_j$, whose coefficients are smooth functions with compact support in $\mb{C}^{n}$, we have
\begin{align*}
(\p\overline{\p}u)(\a)=\int_{\mb{C}^{n}}u
\sum^n_{i,j=1}\frac{\p^2\a_{i\overline{j}}}{\p z_i\p\overline{z}_j}.
\end{align*}
Let $\overline{\p}u=v=\sum^n_{j=1}v_jd\overline{z}_j\in
L^2_{(0,1)}(\mb{C}^{n},e^{-\vp})$. For any $(0,1)$ test form $\b=\sum^n_{j=1}\b_jd\overline{z}_j$, whose coefficients are smooth functions with compact support in $\mb{C}^{n}$, we have
$$-\int_{\mb{C}^{n}}u\sum^n_{j=1}
\frac{\p\b_j}{\p\overline{z}_j}=(\overline{\p}u)(\b)=v(\b)=\int_{\mb{C}^{n}}\sum^n_{j=1}v_j\b_j.$$
Then
\begin{align*}
(\p(\overline{\p}u))(\a)=(\p v)(\a)=-\int_{\mb{C}^{n}}
\sum^n_{i,j=1}v_j\frac{\p^2\a_{i\overline{j}}}{\p z_i}
=-\int_{\mb{C}^{n}}
\sum^n_{j=1}v_j\lt(\sum^n_{i=1}\frac{\p^2\a_{i\overline{j}}}{\p z_i}\rt)=\int_{\mb{C}^{n}}u
\sum^n_{i,j=1}\frac{\p^2\a_{i\overline{j}}}{\p z_i\p\overline{z}_j}
\end{align*}
if let $\b_j=\lt(\sum^n_{i=1}\frac{\p^2\a_{i\overline{j}}}{\p z_i}\rt)$. Thus, $\p\overline{\p}u=\p(\overline{\p}u)$.

Using the same method, we can prove that if $\p u\in L^2_{{1,0}}(\mb{C}^{n},e^{-\vp})$, then $\overline{\p}\p u=\overline{\p}(\p u)$.
\ep

To prove the main theorem, we also need the following simple version of H\"{o}rmander Theorem \cite{4} (page 92, Lemma 4.4.1 with $\vp=|z|^2$).

\bigskip
\noindent\textbf{H\"{o}rmander Theorem.} (A simple version for $(0,1)$ forms)
\textit{For each $f\in L^2_{(0,1)}(\mb{C}^n,e^{-|z|^2})$ such that $\overline{\p}f=0$, there exists a solution $u$ in $L^2(\mb{C}^n,e^{-|z|^2})$ solving equation
$$\overline{\p}u=f$$ in $\mb{C}^n$, in the sense of distributions,
with the norm estimate $$\int_{\mb{C}^n} |u|^2e^{-|z|^2}\leq
2\int_{\mb{C}^n}|f|^2e^{-|z|^2}.$$}

Now we are ready to give the proof of the main theorem.

\bp
For $f\in L^2_{(1,1)}(\mb{C}^n,e^{-|z|^2})$, by Lemma \ref{wj1} we have
\begin{align}\label{23}
f=f_1+\sqrt{-1}f_2,
\end{align}
where $f_1,f_2\in L^2_{2}(\mb{R}^{2n},e^{-|x|^2})$. For $f_1$, by Poincar\'e Lemma ($\vp=|x|^2$, $c=2$, $p=1$) on $\mb{R}^{2n}$, there exists $v_1\in L^2_{1}(\mb{R}^{2n},e^{-|x|^2})$ such that
\begin{align}\label{15}
dv_1=f_1
\end{align}
with
\begin{align}\label{zygs1}
\int_{\mb{R}^{2n}} |v_1|^2e^{-|x|^2}\leq
\frac{1}{4}\int_{\mb{R}^{2n}}|f_1|^2e^{-|x|^2}.
\end{align}
For $v_1$, by Lemma \ref{wj2} we have
\begin{align}\label{16}
v_1=v_1^{1,0}+v_1^{0,1},
\end{align}
where $v_1^{1,0}\in L^2_{{1,0}}(\mb{C}^{n},e^{-|z|^2})$, $v_1^{0,1}\in L^2_{{0,1}}(\mb{C}^{n},e^{-|z|^2})$, $\overline{v_1^{1,0}}=v_1^{0,1}$ and $\overline{v_1^{0,1}}=v_1^{1,0}$.
By (\ref{15}) and (\ref{16}), we have
\begin{align}\label{17}
f_1
=(\p+\overline{\p})(v_1^{1,0}+v_1^{0,1})
=\p v_1^{1,0}+\p v_1^{0,1}+\overline{\p}v_1^{1,0} +\overline{\p}v_1^{0,1}.
\end{align}
Note that $\p v_1^{1,0}$ is a $(2,0)$ form, $\overline{\p}v_1^{0,1}$ is a $(0,2)$ form and $f_1=\frac{1}{2}(f+\overline{f})$ can be seen as a $(1,1)$ form.
So from (\ref{17}), we have $\p v_1^{1,0}=0$, $\overline{\p}v_1^{0,1}=0$ and
\begin{align}\label{18}
\p v_1^{0,1}+\overline{\p}v_1^{1,0}=f_1.
\end{align}
For $v_1^{0,1}$, by H\"{o}rmander Theorem, there exists $u_1\in L^2(\mb{C}^{n},e^{-|z|^2})$ such that
\begin{align}\label{19}
\overline{\p}u_1=v_1^{0,1},
\end{align}
with
\begin{align}\label{zygs2}
\int_{\mb{C}^{n}} |u_1|^2e^{-|z|^2}\leq
2\int_{\mb{C}^{n}}\lt|v_1^{0,1}\rt|^2e^{-|z|^2}.
\end{align}
So for $u_1$, by Lemma \ref{pu} and $\overline{v_1^{0,1}}=v_1^{1,0}$, we have
\begin{align}\label{20}
\p\overline{u_1}=\overline{\overline{\p}u_1}=v_1^{1,0}.
\end{align}
Then by Lemma \ref{wj2}, \ref{huhuan}, \ref{zuijiao}, (\ref{zygs1}), and (\ref{18})-(\ref{20}), we obtain
\begin{align}\label{21}
\p\overline{\p}\lt(u_1-\overline{u_1}\rt)=
\p\overline{\p}u_1-\p\overline{\p}\overline{u_1}=
\p\overline{\p}u_1+\overline{\p}\p\overline{u_1}
=\p(\overline{\p}u_1)+\overline{\p}(\p\overline{u_1})=\p v_1^{0,1}+\overline{\p}v_1^{1,0}=f_1,
\end{align}
with
\begin{align}
\int_{\mb{C}^{n}} |u_1-\overline{u_1}|^2e^{-|z|^2}&\leq
4\int_{\mb{C}^{n}}|u_1|^2e^{-|z|^2}
\leq8\int_{\mb{C}^{n}}\lt|v^{0,1}_1\rt|^2e^{-|z|^2}\nonumber\\
&=2\int_{\mb{R}^{2n}} |v_1|^2e^{-|x|^2}\leq
\frac{1}{2}\int_{\mb{R}^{2n}}|f_1|^2e^{-|x|^2}.\label{zygs3}
\end{align}
By the same method for $u_1$, we can prove that there exists $u_2\in L^2(\mb{C}^{n},e^{-|z|^2})$ such that
\begin{align}\label{22}
\p\overline{\p}\lt(u_2-\overline{u_2}\rt)=f_2.
\end{align}
with
\begin{align}\label{zygs3}
\int_{\mb{C}^{n}} |u_2-\overline{u_2}|^2e^{-|z|^2}\leq
\frac{1}{2}\int_{\mb{R}^{2n}}|f_2|^2e^{-|x|^2}.
\end{align}
Let
\begin{align}\label{24}
u=\lt(u_1-\overline{u_1}\rt)+\sqrt{-1}\lt(u_2-\overline{u_2}\rt).
\end{align}
Then by Lemma \ref{wj1}, (\ref{23}) and (\ref{21})-(\ref{24}), we have $u\in L^2(\mb{C}^{n},e^{-|z|^2})$ such that
\begin{align*}
\p\overline{\p}u
=\p\overline{\p}\lt(u_1-\overline{u_1}\rt)
+\sqrt{-1}\p\overline{\p}\lt(u_2-\overline{u_2}\rt)
=f_1+\sqrt{-1}f_2=f,
\end{align*}
with
\begin{align*}
\int_{\mb{C}^{n}} |u|^2e^{-|z|^2}&=\int_{\mb{C}^{n}} |u_1-\overline{u_1}|^2e^{-|z|^2}+\int_{\mb{C}^{n}} |u_2-\overline{u_2}|^2e^{-|z|^2}\\
&\leq\frac{1}{2}\int_{\mb{R}^{2n}}\lt(|f_1|^2+|f_2|^2\rt)e^{-|x|^2}
\\
&=
2\int_{\mb{C}^{n}}|f|^2e^{-|z|^2}.
\end{align*}
\ep

\section{Why Poincar\'e Lemma}
In this section, we explain why we have to use Poincar\'e Lemma in the proof of the main theorem. Naturely, we could have studied the operator $\p\overline{\p}$ in the following sequence of Hilbert spaces
$$ L^2(\mb{C}^n,e^{-\vp})\stackrel{T}{\longrightarrow} L^2_{(1,1)}(\mb{C}^n,e^{-\vp})\stackrel{S}{\longrightarrow} L^2_{(2,2)}(\mb{C}^n,e^{-\vp}),$$
where operators $T$ and $S$ are extensions of $\p\overline{\p}$ in terms of distributions with domains $Dom(T)$ and $Dom(S)$.
Then we could consider the Hilbert space adjoint $T^*$ and then prove the following formula (whose lengthy calculation is omitted).

\bl\label{lemmaH1}
Let $\vp=|z|^2$.
For any smooth $(1,1)$ form $\a$ with compact support in $\mb{C}^n$, we have
\begin{align*}
\lt\|T^*\a\rt\|^2_\vp
&=\lt\|\a\rt\|_\vp^2
+\lt\|\p\overline{\p}\a\rt\|^2_\vp-\lt\|\p\a\rt\|^2_\vp
-\lt\|\overline{\p}\a\rt\|^2_\vp\\
&\ \ \ \ -\sum^n_{i,j,k,l=1}\lt\|\frac{\p^2\a_{i\overline{j}}}{\p z_k\p \overline{z}_l}\rt\|_\vp^2
+\sum^n_{i,j,k,l=1}\int_{\mb{C}^n}\frac{\p^2\a_{i\overline{j}}}{\p z_k\p \overline{z}_l}\overline{\lt(\frac{\p^2\a_{i\overline{l}}}{\p z_k\p \overline{z}_j}
+\frac{\p^2\a_{k\overline{j}}}{\p z_i\p \overline{z}_l}
\rt)}e^{-\vp}\\
&\ \ \ \
+\sum^n_{i,k,l=1}\lt\|\frac{\p\a_{i\overline{l}}}
{\p z_k}\rt\|^2_\vp
+\sum^n_{j,k,l=1}\lt\|\frac{\p\a_{k\overline{j}}}{\p \overline{z}_l}\rt\|^2_\vp.
\end{align*}
\el

Using the same argument as for Lemma \ref{bijin}, we can prove that the set of smooth $(1,1)$ forms with compact support is dense in $Dom(T^*)\cap Dom(S)$ in the graph norm
$$\|\a\|_{\vp}+\|T^*\a\|_{\vp}+\|\p\overline{\p}\a\|_{\vp}.$$
To apply the density argument, we will have to run into a difficulty from using Lemma \ref{bijin}. For example, the term $-\sum^n_{i,j,k,l=1}\lt\|\frac{\p^2\a_{i\overline{j}}}{\p z_k\p \overline{z}_l}\rt\|_{\vp}^2$ is nonpositive and could not be thrown away
before the limiting argument. So if we take limit from the density
of compact support forms, we will have to end up with the square of a distribution, which is obviously absurd in general.

\end{document}